\documentclass[12pt]{amsart} 
\usepackage{amscd, amssymb, epic, eepic}
 
\theoremstyle{plain}
\newtheorem{proposition}{Proposition}

\newtheorem*{theorem*}{Theorem}
\newtheorem{corollary}{Corollary}

\newtheorem*{main}{Main~Theorem}

\theoremstyle{definition}

\theoremstyle{remark}
\newtheorem*{remark}{Remark}

\newtheorem*{acknowledgment}{Acknowledgments}

\newcommand{\sk}{{\mathcal{K}}}
\newcommand{\pp}{\mathbb{P}}
\newcommand{\zz}{\mathbb{Z}}

\newcommand{\cc}{\mathbb{C}}

\newcommand{\mbar}{{\overline{M}}}
\newcommand{\so}{{\mathcal{O}}}

\newcommand{\lp}{\left(}
\newcommand{\rp}{\right)}

\newcommand{\spec}{{\operatorname{Spec}}}

\newcommand{\tr}{\operatorname{tr}}

\newcommand{\ch}{\operatorname{ch}}
\newcommand{\Td}{\operatorname{Td}}
\newcommand{\h}{\mathcal{H}}
\newcommand{\tlb}{\underline{1}}

\begin{document}
\title{Orbifold Euler characteristics of universal cotangent line bundles on 
$\mbar_{1,n}$}
\author{Yuan-Pin Lee}
\address{Department of Mathematics \\
	University of Utah\\
	 Salt Lake City, Utah  84112-0090}
\email{yplee@math.utah.edu}

\thanks{Research partially supported by NSF grants.}


\pagestyle{plain}

\begin{abstract}
A scheme of computing $\chi(\mbar_{1,n}, L_1^{\otimes d_1}\otimes
\cdots \otimes L_n^{\otimes d_n})$ is given. Here $\mbar_{1,n}$ is the 
moduli space of $n$-pointed stable curves of genus one and $L_i$ are the
universal cotangent line bundles defined by $x_i^*(\omega_{C/M})$, where
$C \to \mbar_{1,n}$ is the universal curve, $\omega_{C/M}$ the relative 
dualizing sheaf and $x_i$ the marked point. 
This work is a sequel to \cite{YL1}.
\end{abstract}

\maketitle

\section{Introduction}
This is a sequel of our previous work \cite{YL1}. The entire program is to 
carry out a $K$-theoretic version of Witten--Kontsevich theory on 
two-dimensional gravity \cite{EW1} \cite{MK}, namely to develop $K$-theory 
of some tautological bundles on the moduli spaces of stable curves. 
In the same time it should serve as examples of quantum $K$-theory
\cite{YL2} as the orbifold Euler characteristics
\begin{equation} \label{e:1}
  \chi(\mbar_{g,n}, L_1^{\otimes d_1}\otimes \cdots \otimes L_n^{\otimes d_n})
\end{equation}
are the quantum $K$-invariants of a point (as target space) 
with gravitational descendents.

More precisely, the Witten--Kontsevich theory of two dimensional gravity
states that a generating function of the intersection indices on moduli
spaces $\mbar_{g,n}$
\[
  \int_{\mbar_{g,n}} c_1(L_1)^{d_1} \cdots c_1(L_n)^{d_n}
\]
is the $\tau$-function for the $KdV$ hierarchy. 
Analogously, one could hope that 
a similar generating function of quantum $K$-invariants \eqref{e:1} could be 
a $\tau$-function for a discrete $KdV$ hierarchy \cite{FR} \cite{DG}.
While the discrete $KdV$ has a ``moduli'' parameter space, 
there is only one quantum $K$-theory. It therefore requires some initial
data to pick the right parameter. The low genera computations in genus zero
\cite{YL1} and the current computation in genus one should serve as a 
starting point. Note that the phenomenon of replacing differential equations 
in quantum cohomology by difference equations in quantum $K$-theory was 
also observed in other examples in quantum K-theory 
\cite{AG} \cite{GL} \cite{YL2} \cite{YL3}.

Let us briefly recount the content of this paper.
Our earlier result \cite{YL1} (see also \cite{RP}) showed that the genus zero
quantum $K$-invariants of a point 
(or Euler characteristics of the line bundles 
$\otimes_{i=1}^n L_i^{\otimes d_i}$) can be written in a closed formula. 
In this article we will give a scheme of computing quantum $K$-invariants
in genus one. 
Unfortunately we are not able to find a closed formula for the generating 
functions. It probably requires some additional insights from theory of 
integrable systems and combinatorics in order to do that.

This computational scheme works in the following steps:

\textbf{First step:} Reduce the (generating functions of) $n$-point 
invariants on $\mbar_{1,n}$
\begin{equation} \label{E:e-0}
  \chi \lp \mbar_{1,n}, \prod_{i=1}^n \frac{1}{1-q_i L_i} \rp 
\end{equation}
to $(n-1)$-points invariants on $\mbar_{1,n}$ (Proposition~\ref{P:e3-3}), 
i.e.~invariants involving at most $n-1$ $L_i$'s.

\textbf{Second step:} Transform the $(n-1)$-point (or less) invariants on 
$\mbar_{1,n}$ to $(n-1)$-point invariants on $\mbar_{1,n-1}$ with similar 
expression (Proposition~\ref{P:e3-1}).

\textbf{Third step:} (initial data)
One can apply the above two procedures until finally
$n=1$. We provide a complete information about the invariants on this 
space and therefore complete this scheme (Proposition~\ref{P:e3-2}).

In fact, we have done more. We have included also tensor powers of the 
\emph{Hodge bundle} 
\[
 \h := H^0(\mbar_{1,n}, \omega)
\]
in our generating invariant \eqref{E:e-0}. Here $\omega$ is the
relative dualizing sheaf of the universal curve over $\mbar_{1,n}$.

In the next section we will start with second step and initial data.
The first step requires a version of Riemann--Roch theorem on orbifolds 
by T.~Kawasaki \cite{TK}. For readers's convenience, we include a 
general discussion of the (fine) moduli space $\mbar_{1,n}$ as an orbifold
and Kawasaki's formula.

For the background as well as the standard definitions and notations used 
here, the readers are referred to \cite{FK} \cite{YL1}.

Note that the techniques in this paper can also be used to compute the
orbifold cohomology of $\mbar_{1,n}$. That will be discussed in a future work.

\begin{acknowledgment} 
I wish to thank A.~Givental and R.~Pandharipande for stimulating discussions.
\end{acknowledgment}

\section{Generalities on orbifolds}
In the previous work \cite{YL1} we have shown how to calculate the genus 
zero correlators, i.e.~$\chi(\mbar_{0,n}, \otimes_{i=1}^n 
L_i^{\otimes d_i})$. In that case the moduli space $\mbar_{0,n}$ is a smooth 
scheme and the techniques employed there are standard techniques in algebraic
geometry. For example, the Hirzebruch--Riemann--Roch theorem was used 
in \cite{YL1}. One might try to apply a similar technique to genus one 
as the next step toward a $K$-theoretic version of Witten--Kontsevich theory.
The first difficulty one would encounter here is that
$\mbar_{1,n}$ is no longer a smooth 
scheme and Hirzebruch--Riemann--Roch does not hold.
This means that the dimensional arguments applied in genus zero case
would not work in the present case. However, $\mbar_{1,n}$ as a moduli 
stack is a smooth Deligne--Mumford stack. 
It is well-known that the smooth stack $\mbar_{1,n}$  
can be represented by an analytic orbifold
\footnote{The term ``orbifold'' 
is sometimes reserved for the space which locally is an \emph{effective}
quotient of a complex space by a finite group. Here however we use this term
for more general spaces, sometimes called \emph{smooth orbispaces} or
\emph{smooth stacks}. In particular our orbifolds do not require 
the local group action to be effective.}
in the complex analytic category (see e.g.~\cite{YR}).  
We will denote them by the the same notation. 
For orbifolds there is a version of Riemann--Roch formula proved by 
T.~Kawasaki \cite{TK}. The idea of Kawasaki is to ``globalize'' 
the equivariant holomorphic Lefschetz fixed point theorem. 
Let us begin our discussion by recalling the classical
holomorphic Lefschetz fixed point theorem.

Let $M$ be a compact complex manifold, $V$ a vector bundle over $X$, and
$G$ a finite group acting on $V \to X$ by automorphism. For any $g \in G$,
let $X^g$ denote the fixed point set of $g$ in $X$, and $(N^*)^g$ be the
conormal bundle of $X^g \subset X$. Then we have
\begin{equation} \label{E:e1-1}
 \sum_p (-1)^p \tr (g : H^p(X,V)) = \frac{1}{|G|} \sum_{g\in G} 
 \chi \lp X^g, \frac{\sum_i \lambda_i V_i}{\prod_j (1- \tau_j (N^*)^g_j)} \rp,
\end{equation}
where $V=\oplus_i V_i$ is a decomposition by the $g$ action and 
$\lambda_i$ the character of $g$ on $V_i$. Similar notations for $(N^*)^g$. 

Now suppose that our orbifold $M$ is a global quotient of an analytic manifold
$\tilde{M}$ by a finite group $G$, and the orbibundle $V$ on $M$ is a vector 
bundle $\tilde{V}$ on $\tilde{M}$. Then the (local) sections of $V$ on
$M$ are the $G$-invariant sections of $\tilde{V}$ on $\tilde{M}$ 
by definition.
In this case the Euler characteristic of $V$ on $M$ will be
\[
  \frac{1}{|G|} \sum_{g \in G} \sum_p (-1)^p \tr(g: H^p(\tilde{M}, \tilde{V})).
\]
We can then apply holomorphic Lefschetz fixed point formula \eqref{E:e1-1}
and get
\begin{equation} \label{E:e1-3}
 \begin{aligned}
  \chi(M,V) = &\frac{1}{|G|} \sum_g \chi \lp \tilde{M}^g, 
   \frac{\sum_i \lambda_i \tilde{V}_i} {\prod_j (1- \tau_j (N^*)^g_j)} \rp \\
  &=\int_{\tilde{M}^g} \frac{\sum_i \lambda_i \ch(\tilde{V}_i) 
	\Td(T\tilde{M}^g)} {\prod_j 1- \tau_j \ch((N_j^*)^g)},
 \end{aligned}
\end{equation} 
where $\sum_i \lambda_i \ch(\tilde{V}_i)$ 
is better interpreted as the equivariant
chern characters, and the denominator should be considered as equivariant
``differential form'' as well. We will use the notation
\[
  \mathcal{T}^g(M,V) := \frac{\sum_i \lambda_i \ch(\tilde{V}_i) 
 	\Td(T\tilde{M}^g)} {\prod_j 1- \tau_j \ch((N_j^*)^g)}
\]
to denote this equivariant differential forms.

In general, an orbifold $M$ is locally a finite group quotient of an
analytic manifold. By definition of orbifolds these local charts should glue
together in a nice way. Therefore we should expect to decompose our orbifold
$M$ into a union of strata, such that each stratum has its (quotient) 
singularity type. More precisely, for any point $p \in M$ consider a local 
chart $(\tilde{U}_p, G_p)$ such that locally at $p$, $M$ is represented as
quotient of $\tilde{U}_p$ by $G_p$. Let $(1)=(h_p^1), (h_p^2), \cdots,
(h_p^{n_p})$ be the conjugacy classes of elements in $G_p$. Consider the set
of pairs:
\[   
   \tilde{\Sigma}M := \{(p,(h_p^i)) | i=1, 2, \cdots,n_p\}.
\]
Take one representative $h_p^i$ of $(h_p^i)$. The pair $(p, (h_p^i))$ 
determines exactly one orbit $[\tilde{p}]$ in $\tilde{U}_p^{h_p^i}$ by
the centralizer $Z_{G_p}(h_p^i)$ of $h_p^i$ in $G_p$. Namely  
$\tilde{\Sigma}M$ locally looks like
\[
  \tilde{\Sigma} U_p \cong \coprod_{i=1}^{n_p} \tilde{U}_p^{h_p^i}/
	Z_{G_p}(h_p^i).
\]
Let $\{\tilde{\Sigma}M_i\}$ be the connected components of $\tilde{\Sigma}M$.
This give $\tilde{\Sigma}M$ an anayltic structure.
Assign a number $m_i$ (\emph{multiplicity}) to each $\tilde{\Sigma}M_i$ by
\[  m_i := \left| \ker \lp Z_{G_p}(g) \to \operatorname{Aut}(\tilde{U}_p^g) 
  \rp \right| , 
\]
where $\tilde{U}_p^{g}/Z_{G_p}(g) \subset \tilde{\Sigma}M_i$.
It is proved in \cite{TK} that the formal sum 
$\sum_{g \in G_p} \mathcal{T}^g(U_p, V)$ defines an ``locally equivariant'' 
differential form on $\tilde{\Sigma}M$, and represents a
cohomology class $\mathcal{T}(M,V)$ in $H^*(\tilde{\Sigma}M, \cc)$.
Using this formulation Kawasaki was able to prove the following theorem.

\begin{theorem*} {\rm \cite{TK}}
Let $V$ be an orbibundle over an orbifold $M$. Then
\[ \chi(M, V)=\sum_i \frac{1}{m_i} \int_{\tilde{\Sigma}M_i} \mathcal{T}(M,V).
\]
\end{theorem*}

Strictly speaking, Kawasaki stated his theorem only for \emph{effective}
orbifolds. With a little effor, one sees that the above theorem applies to 
general orbifolds (e.g.~$\mbar_{1,1}$) as well.

In our case, a generic point $p$ in $\mbar_{1,n}, n\ge2$
is smooth, i.~e.~$G_p=1$. Thus there is a component
$\tilde{\Sigma}M_0 \subset \tilde{\Sigma} M$ 
which represents the fundamental class $[\mbar_{1.n}]$.
The contribution from this term in the Riemann-Roch formula will be called
the \emph{main contribution}, and the rest \emph{orbifold contribution}.

\section{Reduction scheme: second and third steps}

All identities of bundles in this section are stated in $K$-theory, though
most hold at the sheaf theoretic level. 
We will use the notations from sheaf theory, e.g.~tensor product etc., when
they really mean the corresponding objects/operations in $K$-theory.
Since we are working on orbifolds (smooth stacks), ``bundles'' here means 
bundles on stacks (orbi-bundles).

Let $\pi : \mbar_{1,n} \to \mbar_{1,n-1}$ be the forgetful map (forgetting
the last marked point). Let $p:\mathcal{C} \to \mbar_{1,n}$ be the universal 
curve and $\omega_p$ the relative dualizing sheaf. Recall that the Hodge 
bundle $\h:=r^0 p_*(\omega_p)$ is a line bundle on $\mbar_{1,n}$. 
The goal of this section is to relate the $n-1$ point invariants
on $\mbar_{1,n}$ to invariants on $\mbar_{1,n-1}$. 

By definition, for any bundle $E$ 
(or more precisely, an element in $K$-theory) 
\[
 \chi(\mbar_{1,n}, E) = P^n_* (E), 
\]
where $P^n : \mbar_{1,n} \to \spec \cc$ is the map to a point, and 
\[
  P^n_*(E) = P^{n-1}_* (\pi_* E)
\]
by the functoriality of push-forwards in $K$-theory. Note that
$P_*$ here stands for \emph{pushforward in $K$-theory}. That is, all higher
direct images are accounted.

It is clear from the above that one only has to compute 
$\pi_* (\otimes L_i^{\otimes k_i})$ in order to achieve the goal.

\begin{proposition} \label{P:e3-1}
\begin{equation} \label{E:e3-1}
 \begin{aligned}
   \ &\pi_* \lp \frac{1}{1-q\h^{-1}}
	\prod_{i=1}^{n-1} \frac{1}{1-q_i L_i} \rp \\
 =  &\frac{1}{1-q \h^{-1}} \Bigg[
   \lp 1-\h^{-1}+\sum_{i=1}^{n-1} \frac{q_i}{1-q_i} \rp \lp \prod_{i=1}^{n-1} 
  \frac{1}{1-q_i L_i} \rp \Bigg] . 
 \end{aligned}
\end{equation}
\end{proposition}

\begin{proof}
First notice that $\h$ is actually equal to $\pi^*(L_1)$, and therefore
the factor $1/(1-q\h)$ commutes with $\pi_*$.
By Riemann--Roch formula, our generating functions 
\eqref{E:e3-1} are rational functions of $q$ and $q_i$'s. Therefore 
\begin{equation} \label{e:6}
  \chi(\frac{1}{1-qL})= \chi(\frac{-q^{-1}L^{-1}}{1-q^{-1}L^{-1}}). 
\end{equation}
Now for $d_i \geq 1$, we have
\begin{align*}
 &\pi_* \lp\otimes_{i=1}^{n-1} L_i^{-d_i}\rp \\
 = &\otimes_{i=1}^{n-1} L_i^{-d_i}
 \lp\tlb - \h^{-1} -\sum_{i,\ d_i \neq 0} 
 \lp\tlb + L_i+ \cdots +L_i^{d_i -1}\rp\rp.
\end{align*}
It is well known (and follows almost immediately from the definition) that 
\begin{equation} \label{e:lemma}
 L_i=\pi^* {L}_i \otimes \so(D_i), 
\end{equation}
on $\mbar_{g,n}$ for any genus. Here the second $L_i$ is the corresponding
line bundle on $\mbar_{1,n-1}$. $D_i$ is the divisor on 
$\mbar_{1,n}$ whose generic elements have the reducible domain curves
of two components, with $x_i,x_{n}$ on one component and the rest on the
other. By projection formula 
\[  \pi_*(\otimes_{i=1}^{n-1} L_i^{-d_i})= \otimes_{i=1}^{n-1} (L_i)^{-d_i} 
   \otimes \pi_* \so(\sum_i -d_i D_i).
\]
From the definition of $K$-theoretic push-forward 
\[ 
  \pi_*(E)= R^0 \pi_*(E) - R^1 \pi_* (E). 
\]
Since $d_i \ge 1$ for all $i$, 
\[
  R^0 \pi_*(\so(\sum_i -d_i D_i))=0. 
\]
Therefore $R^1:= R^1 \pi_* (\so(\sum_i -d_i D_i))$ is a vector bundle. 
We are now reduced to computing $R^1$. 

Now by Serre duality $(R^1)^* = R^0 \pi_* (\omega(\sum_i d_i D_i))$, where 
$\omega$ is the dualizing sheaf. Fiberwisely, $H^0(C, \omega(\sum_i d_i x_i))$ 
is the holomorphic differential with poles of order at most $d_i$ at $x_i$.
Thus we have a filtration of degrees of poles at each marked point $x_i$:
\[
  F_1 \subset F_2 \subset \cdots \subset F_{d_i}
\] 
and the graded bundles
$F_{k+1}/F_k$ is isomorphic to $L_i^{\otimes - k}$ and $F_0=\h$.
Therefore
\[ 
  \begin{split}
  &\pi_*(\otimes_{i=1}^{n-1} L_i^{-d_i}) \\
  =&\otimes_{i=1}^{n-1} (L_i)^{-d_i} \otimes \lp - R^1 \pi_* 
   \so(\sum_i -d_i D_i) \rp  \\
  =&\otimes_{i=1}^{n-1} (L_i)^{-d_i} \otimes \lp - (R^0 \pi_* 
   \so(\sum_i -d_i D_i))^* \rp  \\
  == &\otimes_{i=1}^{n-1} L_i^{-d_i} \lp\tlb - \h^{-1} -\sum_{i,\ d_i \neq 0} 
 \lp\tlb + L_i+ \cdots +L_i^{d_i -1}\rp \rp.
  \end{split}
\]
Notice that the last equality holds only in $K$-theory (using graded objects). 
This means
\begin{align*}
 &\pi_* \lp \prod_{i=1}^{n-1} \frac{q_i^{-1}L_i^{-1}}{1-q_i^{-1} L_i^{-1}}\rp\\
  =&\lp 1-\h^{-1}+\sum_{i=1}^{n-1} \frac{q_i^{-1}}{1-q_i^{-1}} \rp 
	\lp \prod_{i=1}^{n-1} \frac{q_i^{-1}L_i^{-1}}{1-q_i^{-1}L_i^{-1}}\rp ,
\end{align*}
which is equivalent to \eqref{E:e3-1}
\end{proof}

\begin{remark}
This proof does not use orbifold computation and actually works for 
higher genus quantum $K$-theory (of non-point target space) as well \cite{YL2}.
It is called \emph{$K$-theoretic string equation}.
\end{remark}

\section{Reduction scheme: third step}

For the step three, we have the following proposition.

\begin{proposition} \label{P:e3-2}
\begin{subequations}
  \begin{align}
   \chi\lp\mbar_{1,1}, \frac{1}{1-q_1 L_1}\rp = &\frac{1}
	{(1-q_1^4)(1-q_1^6)}, \label{E:e3-2a}\\
   \chi\lp\mbar_{1,1}, \frac{1}{1-q_1 L_1^{-1}}\rp = &\frac{1-q_1^4 -q_1^6}
	{(1-q_1^4)(1-q_1^6)} \label{E:e3-2b}.
  \end{align}
\end{subequations}
\end{proposition}

\begin{proof}
The formula \eqref{E:e3-2a} follows from the theory of modular
forms. The sections of $L_i^{\otimes k}$ are the modular forms of weight 
$2k$. \eqref{E:e3-2a} is a rephrase of the classical result which states that
the space of modular forms are generated by a weight four and a weight six
modular forms. The equation \eqref{E:e3-2b} follows from \eqref{E:e3-2a}
and \eqref{e:6}.

Alternatively, $\chi\lp\mbar_{1,1}, \frac{1}{1-q_1 L_1^{\pm}}\rp$ can be
computed from (generalized version of) Kawasaki's theorem.
\end{proof}

\begin{corollary} \label{C:e3-2}
\[ 
\begin{split}
 &\chi\lp \frac{1}{(1-qL_1^{-1})(1-q_1 L_1)} \rp  \\
 = &\frac{1}{1-q q_1} \lp \frac{1}{(1-q_1^4)(1-q_1^6)} -\frac{q^{10}}
	{(1-q^4)(1-q^6)} \rp .
\end{split}
\]
\end{corollary}

\begin{proof}
Since $1/(1-qL_1^{-1}) = -q^{-1} L_1/(1-q^-1 L_1)$, this formula follows from
Proposition~\ref{P:e3-2}.
\end{proof}

\section{Reduction scheme: first step}
\begin{proposition} \label{P:e3-3}
\begin{align*}
 &\chi\lp \mbar_{1,n}, \frac{1}{1-q \h^{-1}}
	\prod_{i=1}^n \lp\frac{1}{1-q_i L_i} - 
 	\frac{1}{1-q_i \tlb}\rp \rp \\
 =&\lp  \prod_{i=1}^n \frac{q_i}{1-q_i}\rp 
       \frac{(n-1)!}{24 (1-q)\prod_i (1-q_i)}\text{ (main contribution) }\\
	  + &\ \Sigma_n \text{ (orbifold contributions)}.
\end{align*}
The orbifold contributions $\Sigma_n= 0$ for $n \geq 5$.
For $n=2,3,4$: 
\begin{equation} \label{E:e3-7}
 \begin{aligned}
   \Sigma_4 = &\lp\prod_{i=1}^4 \frac{q_i}{1-q_i}\rp 
  \frac{1}{(1+q) \prod_{i=1}^4 (1+q_i)}\\
      &\qquad \qquad \left[  \frac{3}{4} +
  \frac{1}{2} \lp \frac{q}{1+q}-\sum_i \frac{q_i}{1+q_i} \rp \right].
 \end{aligned}
\end{equation}
\begin{subequations} 
 \begin{align}
  \Sigma_3 = &\lp\prod_{i=1}^3 \frac{q_i}{1-q_i}\rp \notag \\
  &\Bigg[ \lp \frac{1}{(1+q) \prod_{i=1}^3 (1+q_i)} \rp
   \lp -\frac{3}{4} + \frac{1}{2}\lp \frac{q}{1+q} -\sum_i  \frac{q_i}{1+q_i}
      \rp \rp \label{E:e3-8a}\\
   + &\frac{1}	{3 (1+q+q^2) \prod_i (1+q_i +q_i^2)} 
     \Big( -1+ q-(2+q)(q_1+q_2+q_3)  \label{E:e3-8b}  \\
   &\qquad -(1+2q)(q_1 q_2+q_2 q_3+q_3 q_1) +(1-q)q_1 q_2 q_3\Big)\notag\Bigg]
 \end{align}
\end{subequations} 
\begin{subequations}
 \begin{align}
  \Sigma_2 = &\lp\prod_{i=1}^2 \frac{q_i}{1-q_i}\rp \notag \\
 \Bigg[ &\lp \frac{1}{(1+q)\prod_i1+q_i}\rp
    \lp \frac{3}{8}  + \frac{1}{4}
  \lp \frac{q}{1+q}-\sum_i \frac{q_i}{1+q_i}\rp\rp \label{E:e3-9a} \\
  + &\frac{1-q+q_1 +q_2 +q q_1+q q_2 -q_1 q_2+ q q_1 q_2}
    {4 (1+q^2)(1+q_1^2) (1+q_2^2)} \label{E:e3-9b} \\
  + &\frac{1-q+2 q_1 +2 q_2 +q q_1+q q_2+q_1 q_2 - q q_1 q_2}
      {3 (1+q+q^2)(1+q_1+q_1^2)(1+q_2+q_2^2)}\Bigg] 
	\label{E:e3-9c}   
 \end{align}
\end{subequations}
\end{proposition}

\begin{remark}
It's straightforward to check that Proposition~\ref{P:e3-3} actually gives 
power series in $q,q_i$ with \emph{integral} coefficients.
\end{remark}

\begin{proof}
\textbf{Main term: }
The \emph{main contribution} comes from the \emph{fundamental class} of
$\mbar_{1,n}$ because the generic point in $\mbar_{1,n}$ has no automorphism
for $n\geq 2$. Notice that the integrand has the form
\begin{equation} \label{E:e3-10}
  (\text{constant})\lp 1-\frac{q}{1-q}c_1(\h)\rp \lp\prod_{i=1}^n c_1(L_i)\rp
   + \text{higher order terms}.
\end{equation} 
But the dimension of $\mbar_{1,n}$ is $n$, which implies that the only term
contributing to the integration is of the form $\prod_i c_1(L_i)$.
Now the main term can be computed
by the \emph{dilation} equation:
\[
 \int_{\left[\mbar_{g,n}\right]} 
	\otimes_{i=1}^{n-1} c_1(L_i^{\otimes d_i}) c_1(L_n)
  = (2(g-1)+n-1)  \int_{\left[ \mbar_{g,n-1}\right]} 
	\otimes_{i=1}^{n-1} c_1(L_i^{\otimes d_i}).
\]

\textbf{Orbifold terms: }
A point in $\mbar_{1,n}$ is a stable curve of genus one, which consists of
one ``elliptic component'' $E_c$ with some (chains of) rational curves 
attached to it. Here $E_c$ could either be a pointed elliptic curve or
a loop formed by pointed rational curves. Therefore the Kawasaki stratum 
will have a form like $\mbar_{E_c} \times \prod \mbar_{0,m_i}$.
The integration over $\prod \mbar_{0,n_i}$ turns out to be zero due to the
smoothness of $\mbar_{0,m}$ and the following dimension counting: 
the integrand \eqref{E:e3-10} has cohomological degree $\ge 2m$, 
the number of special points in a rational component. However, 
$\dim_{\cc}(\mbar_{0,m}) = m-3$.

We conclude from above that the only orbifold contribution $\Sigma_n$ 
come from those Kawasaki stratum which is of the form $\mbar_{E_c}$.
It is easy to see that automorphisms occur only in the following two cases:
either $E_c$ is a smooth elliptic irreducible component 
or $E_c$ has only two rational components, each carrying two marked 
points, as shown in the following graph (call it ($*$)).

\begin{center}
\setlength{\unitlength}{0.00083333in}
\begingroup\makeatletter\ifx\SetFigFont\undefined
\def\x#1#2#3#4#5#6#7\relax{\def\x{#1#2#3#4#5#6}}%
\expandafter\x\fmtname xxxxxx\relax \def\y{splain}%
\ifx\x\y   
\gdef\SetFigFont#1#2#3{%
  \ifnum #1<17\tiny\else \ifnum #1<20\small\else
  \ifnum #1<24\normalsize\else \ifnum #1<29\large\else
  \ifnum #1<34\Large\else \ifnum #1<41\LARGE\else
     \huge\fi\fi\fi\fi\fi\fi
  \csname #3\endcsname}%
\else
\gdef\SetFigFont#1#2#3{\begingroup
  \count@#1\relax \ifnum 25<\count@\count@25\fi
  \def\x{\endgroup\@setsize\SetFigFont{#2pt}}%
  \expandafter\x
    \csname \romannumeral\the\count@ pt\expandafter\endcsname
    \csname @\romannumeral\the\count@ pt\endcsname
  \csname #3\endcsname}%
\fi
\fi\endgroup
{\renewcommand{\dashlinestretch}{30}
\begin{picture}(2341,1713)(0,-10)
\put(1170.500,398.500){\arc{2426.160}{3.6006}{5.8241}}
\put(1170.500,1986.000){\arc{2401.172}{0.2526}{2.8890}}
\path(608,1611)(683,1386)
\path(1658,1686)(1508,1461)
\path(833,936)(758,786)
\path(1658,936)(1733,786)
\put(383,36){\makebox(0,0)[lb]{\smash{{{\SetFigFont{12}{14.4}{rm}$E_c$ 
is of ($*$) type}}}}}
\end{picture}
}
\end{center}

Let us find all possible genus one stable curves with $n \ge 2$ marked points
which has nontrivial automorphisms.  They belongs to one of the following 
3 types.
(All marked points must be positioned on the fixed points of 
$\operatorname{Aut(E_c,x)}$.)

\begin{enumerate}
 \item The generic element is $E_c$, a genus one irreducible curve with 
  at most four marked points at four 2-torsion points ($\zz/2$ symmetry).
  Note that the corresponding strata contain three (special) moduli points of 
  the form ($*$) with marked points placed at fixed points of $\zz/2$ action.
 \item $E_c$ is a genus one irreducible curve with $\tau=e^{\pi/3}$ such that
	(at most) three marked points sit at fixed points of $\zz/3$ symmetry.
 \item $E_c$ is a genus one irreducible curve with $\tau=i$ with 
	(at most) two marked points at fixed points of $\zz/4$ symmetry.
\end{enumerate}

Now it's obvious that if $n \geq 5$, $\Sigma_n =0$.
Let us deal with the cases $n=2,3,4$ separately.

$\mathbf{n=4 :}$ Only type (1) is possible. The only nontrivial contribution
of singular stratum comes from type (1). But the singular stratum itself is 
the moduli space of elliptic curve with all two-torsion points marked, 
i.~e.~$\pp^1$. We can apply Kawasaki's formula:
\begin{align*}
 \Sigma_4 
 =&\prod_{j=1}^4 \frac{q_j}{1-q_j}\int_{[\pp^1]} \frac{\Td(\pp^1)}
  {\prod_{k=1}^3 \lp 1+e^{c_1(N^*_k)}\rp }
	\lp\frac{1}{1+q e^{c_1(\h^{-1})}} \rp 
    \prod_{j=1}^4 \frac{-e^{c_1(L_j)} -1} {1+q_j e^{c_1(L_j)}} \\
 \intertext{$N^*=\oplus_j N^*_j$ is the conormal bundle (of rank three). 
 	$\zz/2$ acts on $N^*$ and $L_j$ by $-1$.}
 =&\prod \frac{q_j}{1-q_j} \int_{[\pp^1]} (1+2 c_1(L))\lp 1+ \frac{1}{2} 
  c_1(\pp^1)\rp \lp 1- \frac{1}{2} c_1(N^*)\rp  \\
  &\qquad \qquad \qquad \cdot \lp\frac{1+\frac{q}{1+q} c_1(\h)}{1+q}\rp
   \prod_{j=1}^4 \frac{1- \frac{q_j}{1+q_j} c_1(L)}{1+q_j} \\
 \intertext{Here $c_1(L)=c_1(L_j)$ for any $i$. Because $c_1(N^*)= 
  c_1(T^*\mbar_{1,4}|_{\pp^1}) - c_1(T^*\pp^1)$, then} 
 =&\prod \frac{q_j}{1-q_j} \int_{[\pp^1]} (1+2 c_1(L)
  - \frac{1}{2} c_1(T^*\mbar_{1,4}|_{\pp^1})) \\
    & \qquad \qquad \qquad \cdot \lp \frac{1+\frac{q}{1+q} c_1(\h)}{1+q} \rp
  \prod_{j=1}^4 \frac{1- \frac{q_j}{1+q_j} c_1(L)}{1+q_j}. 
\end{align*}
We claim that $\int_{\pp^1} c_1(L)= 1/2$, $\int_{\pp^1} c_1(\h)=1/2$, 
$\int_{\pp^1} c_1(T^*\mbar_{1,4})= -1 + 3/2= 1/2$. 
This can be seen by a local computation.
First for the generic (smooth) curves the formal cotangent space is 
represented
by $H^0(K_E^{\otimes 2}(x_1+x_2+x_3+x_4))$, and $L^{\otimes 2} = \so(1)$.
Second, the forgetful morphism from $\mbar_{1,4}$ to $\mbar_{1,1}$ restricted
to this singular strata $\pp^1$ is a $\zz/2$-equivariant morphism. 
Thus $\h=L_1$ on this strata. Third, these bundles
are trivial over the Zariski open set $U$ formed by smooth curves. One
can trivialize these bundles on $U$ and present bundles on this 
strata $\pp^1$ by choosing transition
functions from $U$ to the neighborhoods of the three degenerate points.
A local computation at degenerate curves show that the generic presentation
of these bundles on $U$ carries through to these three degenerate points.
Combining all above we see that $\Sigma_4=\eqref{E:e3-7}$.

$\mathbf{n=3 :}$ Type (1) and type (2) are possible. 
Since three points out of four
two-torsion  points are sufficient to determine a basis of the two-torsion
points in the Jacobian, the strata for type (1) is still $\pp^1$. 
The contribution of type (1) to $\Sigma_3$ is
\begin{align*}
 &\frac{1}{2} \prod_{j=1}^3 \frac{q_j}{1-q_j} 
   \int_{[\pp^1]} (-8) \lp 1+\frac{3}{2} c_1(L)\rp 
   \lp 1+\frac{1}{2} c_1(T \pp^1)\rp  \\
  &\qquad \qquad \cdot \lp \frac{1}{4}\rp \lp 1-\frac{1}{2} c_1(N^*)\rp 
  \frac{1+\frac{q}{1+q} c_1(\h)}{1+q}
   \prod_j \frac{1- \frac{q_j}{1+q_j} c_1(L)} {1+q_j} \\
 =&\lp \prod \frac{q_j}{1-q_j}\rp \frac{1}{\prod_j (1+q_j)} \Bigg[-\frac{3}{4}
   +\frac{1}{2}\lp \frac{q}{1+q}- \sum_j \frac{q_j}{1+q_j}\rp \Bigg]. \\
 =& \eqref{E:e3-8a}
\end{align*}
The fixed point set of type (2) is just one point. The contribution is ($\eta=
e^{2\pi/3}$):
\begin{align*}
 &\frac{1}{3} \prod_{j=1}^3 \frac{q_j}{1-q_j} 
  \Bigg[ \frac{(\eta-1)^3}{(1-\eta^{-4})(1-\eta^{-2})^2 
  (1-q\eta^2) \prod_j (1-q_j \eta)} \\
     &\qquad \qquad
      + \frac{(\eta^2-1)^3}{(1-\eta^{-2})(1-\eta{-1})^2 
     (1-q\eta)\prod_j (1-q_j \eta^2)}\Bigg] \\
 =& \eqref{E:e3-8b}
\end{align*}

$\mathbf{n=2 :}$ Now all three types can contribute. 
For type (1) the fixed point locus is $\pp^1/(\zz/2)$. 
The contribution to Kawasaki's formula is
\begin{align*}
 &\prod_{j=1}^2 \frac{q_j}{1-q_j} \int_{[{\pp^1}/{(\zz/2)}]} 
  \lp 1+c_1(L) -\frac{1}{2} c_1(T^* \mbar_{1,2}|\lp \pp^1/2\rp ) \rp \\
  &\qquad \qquad \cdot \frac{1+ \frac{q}{1+q}c_1(\h)}{1+q}
  \prod_j  \frac{1- \frac{q_j}{1+q_j} c_1(L_j)}{1+q_j} \\
 =&\prod \frac{q_j}{1-q_j^2} \Bigg[ \frac{3}{8} + \frac{1}{4} \lp 
  \frac{q}{1+q} - \sum_j \frac{q_j}
   {1+q_j}\rp \Bigg] \\
 =& \eqref{E:e3-9a}
\end{align*}
The type (2) contribution is: ($\eta=e^{2\pi/3}$)
\begin{align*}
 &\lp \prod \frac{q_j}{1-q_j}\rp \frac{1}{3} 
  \Bigg[\frac{(\eta-1)^2}{(1-\eta^2)(1-\eta) (1-q\eta^2)
  \prod_j (1-q_j \eta)} \\
   &\qquad \qquad +\frac{(\eta^2-1)^2}{(1-\eta)(1-\eta^2) 
  (1-q\eta)\prod_j (1-q_j \eta^2)} \Bigg] \\
 =& \eqref{E:e3-9b}
\end{align*}
The type (3) contribution is ($i=\sqrt{-1}$):
\begin{align*}
 &\lp\prod \frac{q_j}{1-q_j} \rp\frac{1}{4} \Bigg[ \frac{1}{(1+1)(1+i)} 
    \frac{(-i-1)^2}{(1-q_j)(1+q_1 i)(1+q_2 i)} \\
 &\qquad \qquad + \frac{1}{(1+1)(1-i)}
   \frac{(i-1)^2}{(1+q_j)(1-q_1 i)(1-q_2 i)} \Bigg]  \\
 =& \eqref{E:e3-9c}
\end{align*}
\end{proof}
 
Combining the above three propositions we then get a scheme for
computing all Euler characteristics 
\begin{equation} \label{E:e3-11}
 \chi \lp\mbar_{1,n}, \h^{-d}\otimes(\otimes_{i=1}^n L_j^{\otimes d_j} )\rp 
\end{equation}
for $d, d_j \geq 0$. However, Kawasaki--Riemann--Roch formula as shown above
implies that all terms involved are \emph{rational functions} of $q$'s. It is
therefore possible to formally change $q_i$ to $Q_i^{-1}$ and write it as
a rational function for $Q_i$.
More precisely, we have got the formula 
\[
  f(d_1,\cdots,d_n) = 
  \chi (\mbar_{1,n}, \h^{-d}\otimes \prod_i L_i^{\otimes d_i})
\]
for $d, d_j \geq 0$. If one write the RHS in Kawasaki--Riemann--Roch formula,
one finds that the only difference between positive
and negative powers of $\h^{-1}$ and $L_j$'s in the integral formula is 
merely a formal change of signs of $d$ and $d_j$. Thus this
scheme applies to any integral powers of $L_j, \h$. Summarizing:

\begin{main}
The above gives an effective scheme to compute 
\[
 \chi( \mbar_{1,n}, \h^d \otimes (\otimes_{i=1}^n L_j^{\otimes d_i}))
\]
for any $d, d_i \in \zz$.
\end{main}

\begin{remark}
The argument of extending the positive coefficient to negative ones
applied to genus zero cases \cite{YL1} as well. This therefore extends our
genus zero result to $L_i^{-1}$.
\end{remark}

\section*{Appendix: A simple proof of vanishing theorem}
\label{SS:appendix}
The purpose of this appendix is to give a very simple and 
\emph{self-contained} proof of Pandharipande's vanishing theorem in genus zero
\[
 H^i (\mbar_{0,n}, \otimes_{i=1}^n L_i^{d_i}) =0
\]
for $i \ge 1$ and $d_i \ge 0$. 
\footnote{The method presented in this appendix can also be used to compute
$H^0(\mbar_{0,n}, \otimes L_i^{\otimes d_i})$.
It is also hoped that this method can help to produce an
$S_n$-equivariant version of our genus zero formula \cite{YL1}, which
is needed in the quantum $K$-theory computation of general target spaces.}

Let us first note that in the case $n=4$, $\mbar_{0,4} = \pp^1$ and
$L_i \cong \so(1)$. The vanishing theorem obviously holds.

Let $V := \oplus_{i=1}^n L_i$. Choose generic sections $s_i$ of $L_i$ such
that the zero loci of $s_i$ are in general position. Consider the Koszul
complex
\[
(\mathcal{K}): \quad 0 {\to} \mathcal{O} \overset{d}{\to} V 
\overset{d}{\to} \bigwedge^2 V \to \cdots 
\overset{d}{\to} \bigwedge^n V \to 0,
\]
where the differential $d$ is defined as $d := \sum_{i=1}^n s_i \wedge$ with
$s_i$ considered as a section of $V$.
Because of our choice of $s_i$, this complex $(\mathcal{K})$ is exact.

Now consider the double complex $(C^p (\underline{U}, \mathcal{K}^q); \delta,
d)$, where $\underline{U}$ is a covering of $\mbar_{0,n}$, $\sk^q := \bigwedge^q
V$, $C^p$ are the \v{C}ech cochain groups, $\delta$ is the \v{C}ech 
differential. 

\begin{center}
\setlength{\unitlength}{0.0005in}
%

\begingroup\makeatletter\ifx\SetFigFont\undefined
\def\x#1#2#3#4#5#6#7\relax{\def\x{#1#2#3#4#5#6}}%
\expandafter\x\fmtname xxxxxx\relax \def\y{splain}%
\ifx\x\y   
\gdef\SetFigFont#1#2#3{%
  \ifnum #1<17\tiny\else \ifnum #1<20\small\else
  \ifnum #1<24\normalsize\else \ifnum #1<29\large\else
  \ifnum #1<34\Large\else \ifnum #1<41\LARGE\else
     \huge\fi\fi\fi\fi\fi\fi
  \csname #3\endcsname}%
\else
\gdef\SetFigFont#1#2#3{\begingroup
  \count@#1\relax \ifnum 25<\count@\count@25\fi
  \def\x{\endgroup\@setsize\SetFigFont{#2pt}}%
  \expandafter\x
    \csname \romannumeral\the\count@ pt\expandafter\endcsname
    \csname @\romannumeral\the\count@ pt\endcsname
  \csname #3\endcsname}%
\fi
\fi\endgroup
{\renewcommand{\dashlinestretch}{30}
\begin{picture}(5412,4063)(0,-10)
\path(1050,1315)(1050,1840)
\path(1080.000,1720.000)(1050.000,1840.000)(1020.000,1720.000)
\put(1125,1465){\makebox(0,0)[lb]{\smash{{{\SetFigFont{10}{14.4}{rm}$d$}}}}}
\path(2925,1315)(2925,1840)
\path(2955.000,1720.000)(2925.000,1840.000)(2895.000,1720.000)
\put(3000,1465){\makebox(0,0)[lb]{\smash{{{\SetFigFont{10}{14.4}{rm}$d$}}}}}
\path(2850,2590)(2850,3115)
\path(2880.000,2995.000)(2850.000,3115.000)(2820.000,2995.000)
\put(2925,2740){\makebox(0,0)[lb]{\smash{{{\SetFigFont{12}{14.4}{rm}$d$}}}}}
\path(975,2590)(975,3115)
\path(1005.000,2995.000)(975.000,3115.000)(945.000,2995.000)
\put(1050,2740){\makebox(0,0)[lb]{\smash{{{\SetFigFont{10}{14.4}{rm}$d$}}}}}
\path(3825,865)(4425,865)
\path(4305.000,835.000)(4425.000,865.000)(4305.000,895.000)
\put(3900,1015){\makebox(0,0)[lb]{\smash{{{\SetFigFont{10}{14.4}{rm}$\delta$}}}}}
\path(3825,2140)(4425,2140)
\path(4305.000,2110.000)(4425.000,2140.000)(4305.000,2170.000)
\put(3900,2290){\makebox(0,0)[lb]{\smash{{{\SetFigFont{10}{14.4}{rm}$\delta$}}}}}
\path(1800,2140)(2400,2140)
\path(2280.000,2110.000)(2400.000,2140.000)(2280.000,2170.000)
\put(1875,2290){\makebox(0,0)[lb]{\smash{{{\SetFigFont{10}{14.4}{rm}$\delta$}}}}}
\path(1800,865)(2400,865)
\path(2280.000,835.000)(2400.000,865.000)(2280.000,895.000)
\put(1875,1015){\makebox(0,0)[lb]{\smash{{{\SetFigFont{10}{14.4}{rm}$\delta$}}}}}
\path(25,490)(25,4015)
\path(55.000,3895.000)(25.000,4015.000)(-5.000,3895.000)
\path(25,490)(4400,490)
\path(4280.000,460.000)(4400.000,490.000)(4280.000,520.000)
\put(-350,3940){\makebox(0,0)[lb]{\smash{{{\SetFigFont{10}{14.4}{rm}$q$}}}}}
\put(4525,40){\makebox(0,0)[lb]{\smash{{{\SetFigFont{10}{14.4}{rm}$p$}}}}}
\put(325,790){\makebox(0,0)[lb]{\smash{{{\SetFigFont{8}{14.4}{rm}
$C^0(\underline{U}, \so)$}}}}}
\put(2375,790){\makebox(0,0)[lb]{\smash{{{\SetFigFont{8}{14.4}{rm}
$C^1 (\underline{U}, \so)$}}}}}
\put(325,2140){\makebox(0,0)[lb]{\smash{{{\SetFigFont{8}{14.4}{rm}
$C^0 (\underline{U}, \so(V))$}}}}}
\put(2375,2140){\makebox(0,0)[lb]{\smash{{{\SetFigFont{8}{14.4}{rm}
$C^1(\underline{U}, \so(V))$}}}}}
\end{picture}
}
\end{center}

Use two canonical filtrations (by $p$ and $q$ respectively), we obtain
two spectral sequences $'E_r^{p,q}$ and $''E_r^{p,q}$ with $E_2$ terms:
\begin{gather*}
'E_2^{p,q} = H^p (\mbar_{0,n}, \mathcal{H}_d^q (\sk^*)) \\
''E_2^{p,q} =H_d^q (H^p(\mbar_{0,n}, \sk^*)),
\end{gather*}
and these two spectral sequences abut to the same hyper-cohomology
\[
\mathbb{H}^*(\mbar_{0,n}, \sk^*).
\]



Consider the case $H^k(\mbar_{0,n}, L)$, where $L:= \otimes_i 
L_i^{d_i}$ with $d_i \ge 1, \forall i= 1, \cdots, n$. 
Define $L' := \otimes_i L_i^{d_i -1}$. We will use the complex $\sk\otimes L'$ 
(with differential induced from the complex $\sk$) to get an induction. 
Since $\sk \otimes L'$ is still exact,
\[
 'E_2^{p,q} = H^p(\mbar_{0,n}, \underbrace{\mathcal{H}^q(\sk \otimes L')}_{=0} )
 =0.
\]
On the other hand $''E_2^{p,q} = H_d^q (H^p(\mbar_{0,n}, \sk\otimes L'))$. 
Because $H^p(\mbar_{0,n}, \sk^q \otimes L') =0$ for $q<n$ by induction, 
$''E_2^{p,q} =0$ if $p \ne 0$ and $q<n$. 

So $''E_r^{p,q}$ degenerates at $r=2$.
But we have $'E_{r\ge 2} =0$ and $''E$ and $'E$ abut to the same thing,
we conclude that $''E_{r \ge 2} =0$ and therefore
\[
 H^p (\mbar_{0,n}, L) = H^P (\mbar_{0,n}, \bigwedge^n \otimes L') =0, \quad 
\forall p \ne 0.
\]

For the case that some $d_i=0$, say $d_n=0$. 
\begin{equation} \label{e:16}
 \begin{split}
  &H^j (\mbar_{0,n}, \otimes_{i=1}^{n-1} L_i^{d_i}) 
 =H^j(\mbar_{0,n},\pi^*(\otimes_{i=1}^{n-1} L_i^{d_i})\otimes 
 	\so(\sum d_i D_i))\\
 =&H^j \left(\mbar_{0,n-1}, (\otimes_{i=1}^{n-1} L_i^{d_i}) 
	\otimes \Bigl(R^0 \pi_* \otimes \so(\sum d_i D_i) 
	- R_1 \pi_* \otimes \so(\sum d_i D_i) \Bigr) \right),
 \end{split}
\end{equation}
where $\pi: \mbar_{0,n} \to \mbar_{0,n-1}$ as before.
The first equality follows from (5). 
The second equality uses the Leray spectral
sequence argument and the induction hypothesis 
\[
 H^l(\mbar_{0,n-1}, \otimes L_i^{\otimes d_i})=0
\]
for $l \ge 1$. Note that we have used the same symbols $L_i$ for line 
bundles on $\mbar_{0,n}$ and on $\mbar_{0,n-1}$.
It is easy to see that $R^1 \pi_* \so(\sum d_i D_i) =0$ as 
$H^1(C, \so_C(D))=0$ for $C$ rational and degree of $\so_C(D)$ positive.
The $R^0_*$ can be computed in a way similar to the proof of Proposition~1
(see also the proof of Proposition~1 in \cite{YL1})
\[
  R^0 \pi_* \otimes \so(\sum d_i D_i) =
  1+ \sum_{i, d_i \ne 0} \sum_{j=1}^{d_i} L_i^{-j}.
\]
It is now clear that $H^j (\mbar_{0,n}, \otimes_{i=1}^{n-1} L_i^{d_i})$
can be written as a summation of $H^j (\mbar_{0,n-1}, \otimes_i L_i^{d_i})$,
which by induction is zero. 
  
This ends our proof.

\end{document}